\documentclass[12pt]{amsart}
\usepackage{latexsym}

\usepackage{amsmath,amssymb,amsfonts,amscd}

\def\CC{\mathbb{C}}

\def\RR{\mathbb{R}}
\def\ZZ{\mathbb{Z}}
\def\delbar{\overline{\partial}}

\def\cal{\mathcal}
\def\Det{\textrm{Det}}
\def\coker{\textrm{coker }}
\def\oA{\overline{A}}
\def\uA{\underline{A}}
\def\oe{\overline{e}}
\def\ue{\underline{e}}
\def\og{\overline{\gamma}}
\def\ug{\underline{\gamma}}
\def\oL{\overline{L}}
\def\uL{\underline{L}}
\def\os{\overline{s}}
\def\us{\underline{s}}
\def\oS{\overline{S}}

\def\ov{\overline{v}}
\def\uv{\underline{v}}
\def\ow{\overline{w}}
\def\uw{\underline{w}}
\def\ox{\overline{x}}
\def\ux{\underline{x}}

\author{ Fr\'{e}d\'{e}ric Bourgeois and Klaus Mohnke}
\title{Coherent Orientations in Symplectic Field Theory}
\address{Fr\'{e}d\'{e}ric Bourgeois, Centre de Math\'ematiques,
Ecole Polytechnique, 91128 Palaiseau Cedex, France}
\email{fbourgeo@math.polytechnique.fr}
\address{Klaus Mohnke, Institut f\"ur Mathematik, Humboldt-Universit\"at zu
Berlin, Unter den Linden 6, D-10099 Berlin, Germany}
\email{mohnke@mathematik.hu-berlin.de}
\date{}
\newtheorem{theorem}{Theorem}

\newtheorem{proposition}[theorem]{Proposition}

\newtheorem{remark}[theorem]{Remark}
\newtheorem{corollary}[theorem]{Corollary}

\begin{document}
\bibliographystyle{plain}

\begin{abstract}
We study the coherent orientations of the moduli spaces of holomorphic curves
in Symplectic Field Theory, generalizing a construction due to Floer and Hofer.
In particular we examine their behavior at multiple closed Reeb orbits
under change of the asymptotic direction. The orientations are determined
by a certain choice of orientation at each closed Reeb orbit, that is
similar to the orientation of the unstable tangent spaces of critical points in
finite--dimensional Morse theory.
\end{abstract}

\maketitle

\section{Introduction and Main Results}

Contact homology and Symplectic Field Theory were introduced by
Eliashberg, Givental and Hofer (see \cite{Eliashberg:ICM} and
\cite{Eliashberg/Givental/Hofer:SFT}). They can be considered
as a generalization of Floer homology and Gromov--Witten invariants
to obtain invariants for contact manifolds via the study of holomorphic
curves in symplectizations. Our exposition is a contribution to the foundation
of these theories.
We construct a gluing consistent  orientation of the moduli spaces
of punctured holomorphic curves in so-called symplectic cobordisms.
The latter are open symplectic manifolds $(W^{2n},\omega)$ of dimension $2n$
such that outside a compact set there are nowhere vanishing, symplectically
conformal, complete vector fields. In other words, the ends are
symplectomorphic to either $([R,\infty) \times M_+,d(e^t\alpha_+))$ or
$((-\infty,R] \times M_-,d(e^t\alpha_-))$, i.e.~the positive or negative
half of the symplectization of contact manifolds $(M_{\pm},\xi_{\pm})$
with $\xi_{\pm}:=\ker(\alpha_{\pm})$.
The symplectic cobordism is equipped with an almost complex
structure $J$ which is compatible to the symplectic structure $\omega$
and to the contact $1$-forms on the ends. This means that
$\omega(\cdot,J\cdot)$ gives a Riemannian structure on $W$, that
near each end the almost complex structure $J$ is translationally invariant,
$J(\frac{\partial}{\partial t})=R_\alpha$ (where $R_\alpha$ is the Reeb
vector field defined by $i_{R_\alpha}d\alpha = 0, \alpha(R_\alpha)=1$),
and $J$ preserves $\xi$.

We construct orientations on the moduli spaces of holomorphic curves,
satisfying some relations under the gluing operation.
An outline of such a construction was
given in \cite{Eliashberg/Givental/Hofer:SFT}, but uses a different
approach than ours.

Let us assume that the flow of the Reeb vector field $R_\alpha$ at each end of
$M$ is non-degenerate, i.e.~for all its closed orbits (including the multiples)
their Poincar\'{e} return maps on $\xi$ have no eigenvalue $1$.
We denote by ${\cal P}_{\alpha}$ the set of all closed orbits of the Reeb
flow, where we consider multiples of an orbit as different elements.
For each orbit $\gamma \in {\cal P}_{\alpha}$, we fix a trivialization of
$\xi$ along $\gamma$. Then the linearized Reeb flow on $\xi$ along $\gamma$
defines a path in the symplectic group $Symp(2n-2;\RR)$, starting at the
identity and ending at a matrix with all eigenvalues different from $1$.
The Conley--Zehnder index, $\mu(\gamma)\in \ZZ$, is the Maslov index
of this path (see \cite{Robbin/Salamon:maslovforpaths}).


Let us consider a disk $D \subset\CC$ as a complex curve with a
puncture at the origin $0$ and a fixed direction at the puncture, namely
$\theta = 0$ where $(\rho,\theta)$ are the standard polar
coordinates of $\CC$. We say that a smooth map
$u:D \setminus\{0\}\longrightarrow \RR \times M$ is asymptotic to $\gamma$ at
$\pm\infty$ if $u(\rho,\theta)=(u_\RR(\rho,\theta),u_M(\rho,\theta))$,
$\lim_{\rho\to 0}u_\RR(\rho,\theta)=\pm\infty$ and
the uniform limit $\lim_{\rho\to 0} u_M(\rho, \theta)$ exists and
parameterizes $\gamma$ in such a way that $\lim_{\rho\to 0} u_M(\rho,0)= z_{\gamma}$.
We thus may think of the additionally introduced direction at the
puncture as a puncture on the boundary of the compactification
$[0,1]\times S^1$ of $D\setminus\{0\}$ via polar coordinates. The map $u$ thus
extends to the compactification and maps the puncture on the boundary to $z_\gamma$.

Let $(\Sigma,j)$ be a Riemann surface, $j$ its conformal structure.
Let $\ox_1,\ldots \ox_{\os}; \ux_1,\ldots,\ux_{\us}\in \Sigma$ be distinct
points of the surface (which we will call {\em punctures}) and
$\ov_1 \in T_{\ox_1}\Sigma,\ldots,\ov_{\os}\in T_{\ox_{\os}}\Sigma; \uv_1 \in
T_{\ux_1}\Sigma,\ldots,\uv_{\us}\in T_{\ux_{\us}}\Sigma$ be fixed unit vectors
in the tangent spaces at the points, which we call {\em directions}.
We say that a map
$u:\Sigma\setminus\{\ox_1,\ldots,\ox_{\os},\ux_1,\ldots,\ux_{\us}\}
\longrightarrow W$  is asymptotic
to $\og_k \in {\cal P}_{\alpha_{+}}$ in $\ox_k$ at $+\infty$ if there
are complex polar coordinates $(\rho,\theta)$
centered at $\ox_k$ such that the directions $\ov_k$ coincide with
$\theta = 0$ and $u$ is asymptotic to $\og_k$ at
$+\infty$. Analogously, we say that $u$ is asymptotic to
$\ug_l \in {\cal P}_{\alpha_{-}}$ in $\ux_l$ at $-\infty$.

Note that the diffeomorphisms of $\Sigma$ act in an obvious way
on the data described above, giving rise to new conformal structures and new holomorphic
maps with the same asymptotics. We are now ready to define the moduli spaces of
holomorphic curves in a symplectic cobordism $(W,J)$.

Pick a closed Riemann surface $\Sigma$, and
closed orbits $\og_1,\dots,\og_{\os}\in {\cal P}_{\alpha_{+}}$
and $\ug_1,\dots,\ug_{\us}\in{\cal P}_{\alpha_{-}}$ of the flow
of the corresponding Reeb vector fields. Define the moduli space
$$
{\cal M}_{W,J}^{\Sigma}(\og_1,\dots,\og_{\os};\ug_1,\dots,\ug_{\us})
$$
to be the set of conformal structures on $\Sigma$ together with
$\os$ positive and $\us$ negative punctures with directions,
and holomorphic maps into $(W,J)$ asymptotic to
$\og_1,\dots,\og_{\os}$ at $+\infty$ at the positive punctures
and to $\ug_1,\dots,\ug_{\us}$ at $-\infty$ at the negative punctures,
modulo the action by the diffeomorphisms of $\Sigma$.

Two symplectic cobordisms $(W_1,\omega_1)$ and $(W_2,\omega_2)$ can be glued
if the negative end of $(W_1,\omega_1)$ has the form
$((-\infty,-R_0] \times M,d(e^t\alpha))$ and the positive end of
$(W_2,\omega_2)$ looks like $([+R_0,,+\infty) \times M, d(e^t\alpha))$.
The glued cobordism  $W_R$, for $R > R_0$, is obtained after
cutting $W_1$ along $\{ -R \} \times M$, $W_2$ along $\{ +R \} \times M$
and identifying the ends with the identity map. We obtain the glued symplectic
structure $\omega_R$ after multiplying $\omega_1$ by $e^R$ and $\omega_2$ by
$e^{-R}$. If the $(W_i,\omega_i)$ are equipped with compatible almost complex
structures that are identical on the ends we glue, then the glued cobordism
$(W_R,\omega_R)$ is naturally equipped with a compatible almost complex
structure $J_R$. Of course, if at least one of the $(W_i,\omega_i)$ is a
symplectization equipped with an $\RR$-invariant almost complex structure,
the glued structures $(W_R,\omega_R,J_R)$ are independent of $R > R_0$.

We can also glue holomorphic curves in the almost complex manifolds $(W_i,J_i)$.
Consider compact subsets
$$
K_1 \subset {\cal M}_{W_1,J_1}^{\Sigma_1} (\og_1,\dots,\og_{\os};
\ug_1,\dots,\ug_{\us-t},\beta_1,\dots, \beta_t)
$$
and
$$
K_2 \subset {\cal M}_{W_2,J_2}^{\Sigma_2}(\beta_1,\dots,\beta_t,
\og'_{t+1},\dots,\og'_{\os'}; \ug'_1,\dots, \ug'_{\us'})
$$
of a pair of moduli spaces. Let $\Sigma = \Sigma_1 \sharp_t \Sigma_2$ be the
Riemann surface obtained by gluing $\Sigma_1$ and $\Sigma_2$ at the punctures
with asymptotics $\beta_1, \ldots, \beta_t$, identifying their asymptotic directions.
If $(W_2, \omega_2, J_2)$ is a symplectization equipped with an $\RR$-invariant
almost complex structure, we can have $\us > t$ since a holomorphic curve in
$(W_1, J_1)$ can be extended in $(W_2,J_2)$ by a cylinder over a closed Reeb orbit;
otherwise, we must have $\us = t$. There is a similar remark for
$(W_1,\omega_1, J_1)$.

One constructs, for $R$ large enough, a diffeomorphism
$$
\Phi_R : K_1 \times K_2 \longrightarrow {\cal M}_{W_R,J_R}^{\Sigma}
(\og_1,\dots,\og_{\os},\og'_{t+1},\dots,\og'_{\os'};
\ug_1,\dots,\ug_{\us-t},\ug'_1,\dots,\ug'_{\us'}) .
$$
In this paper, we will construct and use a linearized version of this {\em gluing
diffeomorphism}.

Each moduli space ${\cal M}_{W,J}^{\Sigma}$ comes naturally equipped with a real line
bundle, called determinant bundle. Whenever the moduli spaces are cut out transversally
by the Cauchy-Riemann equation (or a small perturbation of these equations), the
determinant bundle is canonically isomorphic to the orientation bundle
$\Lambda^{\rm max}T{\cal M}_{W,J}^{\Sigma}$ of ${\cal M}_{W,J}^{\Sigma}$.
In particular, an orientation of the determinant bundle of ${\cal M}_{W,J}^{\Sigma}$ is
equivalent to an orientation of ${\cal M}_{W,J}^{\Sigma}$. Using gluing analysis of
Fredholm operators, similar to that of Floer and Hofer in \cite{Floer/Hofer:orient},
we will construct orientations for the determinant bundles of the moduli spaces :

\begin{theorem}\label{gluing}
The determinant bundles of the moduli spaces
${\cal M}_{W,J}^{\Sigma}$ are orientable in such a way that the
gluing diffeomorphisms preserve the orientations (up to a sign due to the
reordering of some punctures, see Corollary~\ref{glue}) and the orientation
of a disjoint union of data is induced by the orientation on each
component (see Proposition~\ref{union}).
\end{theorem}

In this article we describe the behavior of the
constructed orientations under the change of certain data.
The first property is a relation for the orientations of the moduli
spaces when reordering  the punctures. This has been known since the
early days of Symplectic Field Theory, for algebraic
reasons : it is necessary for the differential $d$ to be well-defined.

\begin{theorem}\label{orientorderings}
Let
\begin{eqnarray*}
\lefteqn{\underline{\iota}_{l,l+1} : {\cal M}_{W,J}^{\Sigma}
(\og_1,\dots,\og_{\os};\ug_1,\dots,\ug_l,\ug_{l+1},\ldots,\ug_{\us})} \\
&\longrightarrow&
{\cal M}_{W,J}^{\Sigma}
(\og_1,\dots,\og_{\os};\ug_1,\dots,\ug_{l+1},\ug_l,\dots,\ug_{\us})
\end{eqnarray*}
be the map induced by exchanging the $l$th and the $(l+1)$th negative punctures.
There is a natural identification of the determinant bundles of the
corresponding linearizations for $[u,j]$ and $\underline{\iota}_{l,l+1}([u,j])$,
which are simply the same. Then  the orientations given in
Proposition~\ref{gluing}
\begin{enumerate}
\item[(i)] agree, if the product
$(\mu({\ug_l})+(n-3))(\mu({\ug_{l+1}})+(n-3))$ is even;
\item[(ii)] disagree, if it is odd.
\end{enumerate}
A similar statement holds for the change of the ordering at the
positive end.
\end{theorem}

The contribution of the paper is to explain this by the consistency
relation under disjoint union. Next we study the behavior of the
orientation under the change of asymptotic directions
at the punctures which are asymptotic to multiply covered
closed Reeb orbits. We say that a closed Reeb orbit $\gamma_m$ is
{\em bad} if it is the $m$-fold covering of some Reeb orbit $\gamma$
and the difference $\mu(\gamma_m)-\mu(\gamma)$ of their Conley--Zehnder
indices is odd. See e.g.~\cite{Ustilovsky:thesis}
for a discussion. If this happens, then the integer $m$ must be even.
Closed Reeb orbits which are not bad are called {\em good}.

\begin{theorem}\label{orientbadorbits}
Consider the moduli space ${\cal M}^{\Sigma}
(\og_1,\dots,\og_{\os};\ug_1,\dots,\ug_{\us})$,
where  $\og_k$ is the $m$-fold covering of a closed Reeb
orbit. There is a natural $\ZZ_m$--action on it which acts
transitively on the possible direction at $\ox_k$.
This action is orientation preserving if and only if
$\og_k$ is good.
A similar statement holds in the case of a negative puncture.
\end{theorem}

At this point, the role of the chosen points $z_{\gamma}$
and of the asymptotic directions becomes clear. Otherwise,
it would be impossible to fix an orientation for
holomorphic curves with at least one puncture asymptotic to a
multiple Reeb orbit with odd behavior.
The necessity of such a result has been clear since Michael Hutchings
observed that in order to have invariance of Symplectic Field Theory,
certain closed Reeb orbits should be ignored.

We will discuss the Fredholm theory for linearized Cauchy-Riemann
operators, corresponding to holomorphic curves in symplectic cobordisms,
in Section~\ref{setup}. Then in Section~\ref{orientation}
we will generalize the construction of coherent orientations
of Floer and Hofer. In Section~\ref{moduli}, we apply it to the moduli
spaces of holomorphic curves in Symplectic Field Theory. Finally,
in Section~\ref{behavior}, we will discuss the even and odd behavior
under the change of asymptotic directions at multiple closed Reeb orbits
described in Theorem~\ref{orientbadorbits}.

{\em Acknowledgements.} The authors are indebted to Yakov
Eliashberg, Michael Hutchings, Helmut Hofer and Kai Cieliebak for enlightening
discussions and interest in their work. The first author was partially supported
by the Fonds National de la Recherche Scientifique and the European Differential
Geometry Endeavour. The second author was  supported by the
DFG. Both authors also want to thank the Mathematics Department of
Stanford University for its warm hospitality.

\section{Fredholm Theory}\label{setup}

Let $\Sigma$ be a closed Riemann surface with conformal structure $j$,
positive punctures $\ox_k$, $k = 1, \ldots, \os$ and negative punctures $\ux_l$,
$l = 1, \ldots, \us$. On a small disk of radius $e^{-R_0}$ centered on a puncture
$\ox_k$ or $\ux_l$, we will use a local complex coordinate
$z = e^{\mp(s+\mbox{i}\theta)}$ vanishing at $\ox_k$ or $\ux_l$.
Then $(s,\theta) \in [R_0,\infty) \times S^1$ for a positive
puncture and $(s,\theta) \in (-\infty,-R_0] \times S^1$ for a negative puncture;
$(s,\theta)$ are cylindrical coordinates on the punctured disk.
We obtain a compactification $\overline{\Sigma}$ of the Riemann surface
$\Sigma$ by extending the coordinate $s$ to $\pm\infty$ at each puncture.

Let $E$ be a symplectic vector bundle over the closed Riemann surface $\Sigma$
together with symplectic linear identifications $E_{\ox_k} \cong E_{\ux_l} \cong
(\RR^2,\omega_0) \oplus (\RR^{2n-2},\omega_0)$. This gives rise to a symplectic
vector bundle $\overline{E}$ over the compactification $\overline{\Sigma}$ with
fixed trivializations on the boundary components.

Let $\beta : \RR \rightarrow [0,1]$ be a smooth function such that $\beta(s) = 0$
if $s < 0$, $\beta(s) = 1$ if $s > 1$ and $0 \le \beta'(s) \le 2$.
Let $\oe_{k,a}$, $a = 1, 2$ be vectors spanning the first summand of
$E_{\ox_k} \cong \RR^2 \oplus \RR^{2n-2}$.  We similarly define vectors
$\ue_{l,a}$.
We define $\Gamma_{\os,\us} \cong \RR^{2\os + 2\us}$ to be the vector space
generated by the sections $\ow_{k,a}(s) = \beta(|s|-R_0) \oe_{k,a}$ and
$\uw_{l,a}(s) = \beta(|s|-R_0) \ue_{l,a}$.

Let $L^{p,d}_k(E)$ be the Sobolev space of sections $\zeta$ of $E$ that are locally
in $L^p_k$ and such that, near each puncture, $e^{\frac{d}{p}|s|}\zeta(s,\theta)$
is in $L^p_k$ with respect to the cylindrical measure $ds \, d\theta$. The Banach
norm on $L^{p,d}_k(E)$ is defined with respect to the measure $e^{d|s|}dsd\theta$
near the punctures.

Let
\begin{eqnarray*}
{\cal S}:= \!\! &\{& \!\! A:[0,2\pi]\longrightarrow \mbox{Symp}(2n-2;\RR) \mid
1 \notin \mbox{spec}(A(2\pi));\\
& & A(0)=\mbox{Id}, \dot{A}(0)A(0)^{-1}=\dot{A}(2\pi)A(2\pi)^{-1}\ \}
\end{eqnarray*}
be the set of  regular paths in the symplectic group.
We will denote the Maslov index  of a path (see \cite{Robbin/Salamon:maslovforpaths}) of
symplectic matrices $A \in {\cal S}$ by $\mu(A)$.
Note that for a regular closed Reeb orbit $\gamma$ with a fixed symplectic
trivialization of $\xi$, the linearization of the Reeb flow along $\gamma$
determines an element $A_{\gamma} \in {\cal S}$.

Now, for $\oA_k , \uA_l \in{\cal S}$ let ${\cal O}
(\Sigma,E;\oA_1,\ldots,\oA_{\os};\uA_1,\ldots,\uA_{\us})$
be the set of continuous linear operators $L : \Gamma_{\os,\us} \oplus L^{p,d}_k(E)
\longrightarrow L^{p,d}_{k-1}(\Lambda^{0,1}(E))$, with $p > 2$, $k \ge 1$ and
$d > 0$, of the following form :
\begin{enumerate}
\item[(i)] in a trivialization of $E$ in the interior of $\Sigma$, with local
coordinate $z = x + \mbox{i}y$, the operator $L$ looks like
$$
\Big( \frac{\partial}{\partial x} + \mbox{J(z)} \frac{\partial}{\partial y}
+ S(z) \Big) (dx - \mbox{i}dy),
$$
where $J(z)$ and $S(z)$ are locally in $L^p_k$ and $J(z)$ is a compatible complex
structure on $E$;
\item[(ii)] in a trivialization of $\overline{E}$ near a puncture, that extends
the fixed identification at the puncture, with the cylindrical coordinates
$(s,\theta)$, the operator $L$ looks like
$$
\Big(\frac{\partial}{\partial s} +
\mbox{J(s,$\theta$)} \frac{\partial}{\partial \theta} +
S(s,\theta)\Big)(ds-\mbox{i}d\theta),
$$
where $S(s,\theta)$ and $J(s,\theta)$ extend continuously to the compactification
$\overline{\Sigma}$ in such a way that $J(s,\theta) - J(\pm \infty, \theta)$ and
$S(s,\theta) - S(\pm \infty, \theta)$ are in $L^{p,d}_{k-1}$, and $J(s,\theta)$
is a compatible complex structure on $\overline{E}$.
\end{enumerate}

The loop of matrices $S(+\infty,\theta)$ splits in
$E_{\ox_k} \cong \RR^2 \oplus \RR^{2n-2}$ as
$$
S(+\infty,\theta)= \left( \begin{array}{cc}
0 & 0 \\
0 & \oS_k(\theta)
\end{array} \right) ,
$$
where
$$
\oS_k(\theta) := -J(+\infty,\theta)\dot{\oA}_k(\theta)
\oA_k(\theta)^{-1}
$$
is a {\em loop} of $(2n-2)\times (2n-2)$--matrices which are symmetric with
respect to the euclidean structure determined by the symplectic and the
complex structure. There is a similar expression for the loop of matrices
$S(-\infty,\theta)$.

Note that the differential operators we consider have degenerate asymptotics,
because the matrices $S$ vanish on the first summand of
$E_{\ox_k} \cong E_{\ux_l} \cong \RR^2 \oplus \RR^{2n-2}$. Hence, we cannot
directly apply to them the usual Fredholm theory as in \cite{Schwarz:thesis}.

\begin{proposition}  \label{Fred}
Let $L \in {\cal O}(\Sigma,E;\oA_1,\ldots,\oA_{\os};\uA_1,\ldots,\uA_{\us})$.
Then
$$
L : \Gamma_{\os,\us} \oplus L^{p,d}_k(E)  \longrightarrow
L^{p,d}_{k-1}(\Lambda^{0,1}(E))
$$
is a Fredholm operator with index
$$
\mbox{ind}(L) = \sum_{i=1}^{\os}(\mu(\oA_i)-(n-1))
- \sum_{j=1}^{\us}(\mu(\uA_j)+(n-1)) + n\chi(\Sigma) + 2c_1(E).
$$
$\chi(\Sigma)$ denotes the Euler characteristic of $\Sigma$.
\end{proposition}

\begin{proof}
First consider the restriction of $L$ to $L^{p,d}_k(E)$.
Note that the $L^{p,d}_k$ spaces are isomorphic to the $L^p_k$ spaces, via
multiplication by a function that is constant away from the punctures and
given by $e^{\frac{d}{p}|s|}$ near the punctures. Conjugating our restricted
operator with these isomorphisms, we obtain an operator $L'$ from $L^p_k$ to
$L^p_{k-1}$. Moreover, near a puncture $\ox_k$ or $\ux_l$, the operator $L'$
looks like
$$
\Big(\frac{\partial}{\partial s} +
\mbox{J(s,$\theta$)} \frac{\partial}{\partial \theta} +
S(s,\theta) \mp \frac{d}p I \Big)(ds-\mbox{i}d\theta) .
$$
In particular, if the operator $L$ has asymptotics $\oA_k, \uA_l \in  {\cal S}$,
the new operator $L'$ has asymptotics $\oA'_k, \uA'_l$, where
\begin{eqnarray*}
\oA'_k(\theta) = e^{-\mbox{i}\frac{d}p \theta} \left( \begin{array}{cc}
I & 0 \\
0 & \oA_k(\theta)
\end{array} \right)
& \textrm{and} &
\uA'_l(\theta) = e^{+\mbox{i}\frac{d}p \theta} \left( \begin{array}{cc}
I & 0 \\
0 & \uA_l(\theta)
\end{array} \right) .
\end{eqnarray*}
As a consequence of this small perturbation,
the asymptotics are not degenerate anymore, and we obtain
a Fredholm operator (see \cite{Schwarz:thesis}). Its index is given by
$$
n(2-2g-\os - \us) + 2c_1(E) + \sum_{k=1}^{\os}\mu(\oA'_k)
-\sum_{l=1}^{\us} \mu(\uA'_l) .
$$
It is easy to see that $\mu(\oA'_k) = \mu(\oA_k) - 1$ and
$\mu(\uA'_l) = \mu(\uA_l) + 1$.

Let us now consider the summand $\Gamma_{\os,\us}$.
The images under $L$ of these sections have their support near a puncture
and decay exponentially near that puncture,
hence $L$ maps the new elements into $L^{p,d}_{k-1}$. The space $\Gamma_{\os,\us}$
is finite dimensional. Thus the operator $L$ is Fredholm.
Its Fredholm index is the sum of the index of the
operator restricted to $L^{p,d}_k(E)$ and the dimension of the
vector space $\Gamma_{\os,\us}$, which is $2(\os + \us)$.
Hence we obtain the desired formula.
\end{proof}

As in \cite{Floer/Hofer:orient}, the sets
${\cal O}(\Sigma,E;\oA_1,\ldots,\oA_{\os};\uA_1,\ldots,\uA_{\us})$
are topological spaces. Moreover, as in Proposition~8 of
\cite{Floer/Hofer:orient}, these spaces are contractible.
Note that unlike in \cite{Floer/Hofer:orient} we allow the complex structures
to vary on the ends here. Since the space ${\cal J}(\overline{E})$
of compatible complex structures is contractible this does not cause
any trouble similar to that of Theorem~2 of \cite{Floer/Hofer:orient}.
The advantage will be that we can make orientation choices for each
closed Reeb orbit together with its fixed trivialization of $\xi$
without further structure.

By Proposition~\ref{Fred}, the topological spaces
${\cal O}(\Sigma,E;\oA_1,\ldots,\oA_{\os};\uA_1,\ldots,\uA_{\us})$
consist of Fredholm operators. Recall that, for every Fredholm operator $L$,
the determinant space of $L$ is defined by $\Det (L) := (\Lambda^{\rm max} \ker L)
\otimes (\Lambda^{\rm max} \coker L)^*$. Moreover, the determinant spaces of
$L \in {\cal O}$ fit together in a determinant bundle $\Det({\cal O})$ over ${\cal O}$.
Since the spaces ${\cal O}$ are contractible, this determinant bundle is
orientable.
In order to orient these determinant bundles in a coherent fashion, we need to
construct a gluing map to relate orientations on different bundles $\Det({\cal O})$.

Consider two operators
\begin{eqnarray*}
K &\in& {\cal O}(\Sigma,E;\oA_1, \ldots, \oA_{\os};
\uA_1,\ldots,\uA_{\us}) ,\\
L &\in& {\cal O}(\Sigma',E';\oA'_1, \ldots,\oA'_{\os'};
\uA'_1,\ldots,\uA'_{\us'}) ,
\end{eqnarray*}
such that $\uA_{\us+1 - m}$ coincides with $\oA'_m$ for $m = 1, \ldots, t
\le \min(\us,\os')$. We want to construct a glued operator
$$
M_R \in {\cal O}(\Sigma'', E'';
\oA_1, \ldots, \oA_{\os},\oA'_{t+1}, \ldots, \oA'_{\os'};
\uA_1,\ldots,\uA_{\us-t},\uA'_1,\ldots,\uA'_{\us'}) .
$$

Near the punctures $\ux_{\us +1 - m}$, pick cylindrical coordinates
$(s_m, \theta_m) \in (-\infty,-R_0] \times S^1$; near the punctures
$\ox'_m$, pick cylindrical coordinates $(s'_m, \theta'_m) \in
[R_0,+\infty) \times S^1$. Choose $R$ much larger than $R_0$. Cut out
the small disks $\{ s_m < -R \}$ from $\Sigma$ and $\{ s'_m > +R \}$ from
$\Sigma'$ and identify their boundaries for $m = 1, \ldots, t$ so that
$\theta_m = \theta'_m$.  We denote
the resulting Riemann surface by $\Sigma''$. Gluing the bundles $E$ and $E'$ using
their trivializations near the punctures, we obtain a bundle $E''$ over
$\Sigma''$. In the neighborhood of the gluing, we obtain cylindrical
coordinates $(s''_m,\theta''_m)$ defined by $s''_m = s_m + R$,
$\theta''_m = \theta_m$ on $\Sigma$ and $s''_m = s'_m - R$,
$\theta''_m = \theta'_m$ on $\Sigma'$.

In cylindrical coordinates $(s''_m,\theta''_m)$, we define
$$
M_R := \Big(\frac{\partial}{\partial s''_m} +
J''(s''_m,\theta''_m) \frac{\partial}{\partial \theta''_m} +
S''(s''_m,\theta''_m)\Big)(ds''_m-\mbox{i}d\theta''_m),
$$
where
\begin{eqnarray*}
S''(s''_m,\theta''_m) = S(-\infty, \theta''_m) & \!\! + & \!\!
\beta(s''_m) (S(s_m,\theta_m)-S(-\infty,\theta_m)) \\
& \!\! + & \!\! \beta(-s''_m) (S'(s'_m,\theta'_m)-S'(+\infty,\theta'_m))
\end{eqnarray*}
and $J''(s''_m,\theta''_m)$, $s''_m \in [-1,+1]$, is a path of compatible
complex structures on $E''$ interpolating between $J(-R+1,\theta_m)$ and
$J'(+R-1,\theta'_m)$.

Note that, on $L_k^{p,d}(E'')$, we will not use the usual Banach norm,
but a modified norm. For $\zeta \in L_k^{p,d}(E'')$, let
$\bar \zeta_m = \frac1{2\pi}\int_{s''_m = 0} \pi_m \zeta \, d\theta''_m$,
where $\pi_m$ is the projection to the first summand of
$E_{\ux_{\us - t + m}} \cong E_{\ox'_m} \cong \RR^2 \oplus \RR^{2n-2}$.
Let ${\bar \beta}(s) = \beta(-s + R - R_0) \beta(s + R - R_0)$.
We define
$$
\| \zeta \|_{L_k^{p,d}(E'')} = \| \zeta - \sum_{m=1}^t {\bar \beta}(s''_m)
{\bar \zeta_m} \|_R + \sum_{m=1}^t | \bar \zeta_m | ,
$$
where $\| \cdot \|_R$ is the Banach norm with respect to a measure of the same
form as before near the punctures but with a measure
$e^{d(R-R_0 - |s''_m|)}ds''_m d\theta''_m$ in the neighborhood of the gluing.

On the other hand, on $L_{k-1}^{p,d}(\Lambda^{0,1}(E''))$,
we will use the Banach norm $\| \cdot \|_R$, with a measure of the same form as
before near the punctures but with a measure
$e^{d(R-R_0 - |s''_m|)}ds''_m d\theta''_m$ in the neighborhood of the gluing.

Let $\Gamma_t$ be the vector space generated by the pairs
$(\uw_{\us+1-m,a}, \ow'_{m,a}) \in \Gamma_{\os,\us} \oplus \Gamma'_{\os',\us'}$,
for $m = 1, \ldots, t$ and $a = 1, 2$. We write $\os'' = \os + \os' - t$ and
$\us'' = \us + \us' - t$. Let $\Gamma''_{\os'',\us''}$ be the vector space
generated by the sections $\ow_{k,a}$, for $k = 1, \ldots, \os $, the sections
$\uw_{l,a}$ for $l = 1, \ldots, \us -t$, the sections $\ow'_{k,a}$, for
$k = t+1, \ldots, \os'$, and the sections $\uw'_{l,a}$, for $l = 1, \ldots, \us'$
and $a=1,2$.
We denote by $K \oplus_{\Gamma_t} L$ the restriction of $K \oplus L$ to
$\Gamma''_{\os'',\us''} \oplus \Gamma_t \oplus L_k^{p,d}(E)
\oplus L_k^{p,d}(E')$.

\begin{proposition} \label{rinv}
Assume that the operator $K \oplus_{\Gamma_t} L$ is surjective.
Then the operator
$$
M_R : \Gamma''_{\os'',\us''} \oplus
L_k^{p,d}(E'') \rightarrow L_{k-1}^{p,d}(\Lambda^{0,1}(E''))
$$
has a uniformly bounded right inverse $Q_R$, if $R$ is sufficiently large.
\end{proposition}

\noindent
In order to prove this proposition, we adapt the gluing construction of
McDuff and Salamon \cite{McDuff/Salamon:Quantum}.

Let $\gamma_R : \RR \rightarrow [0,1]$ be a smooth decreasing function such that
$\gamma_R(s) = 0$ for $s \ge \frac{R-R_0}2$, $\gamma_R(s) = 1$ for
$s \le 1$ and all derivatives of $\gamma_R$ uniformly converge to $0$ as
$R \rightarrow \infty$.

Let us define the gluing map
\begin{eqnarray*}
g_R : \Gamma_t \oplus L_k^{p,d}(E) \oplus L_k^{p,d}(E')
&\rightarrow& L_k^{p,d}(E'') \\
(v, \zeta ,\zeta') &\rightarrow& \zeta'' = \zeta^0 + v
\end{eqnarray*}
where, in the neck with coordinates $(s''_m,\theta''_m)$,
$$
\zeta^0(s''_m,\theta''_m) = \left\{ \begin{array}{ll}
\zeta(s''_m,\theta''_m) +
\gamma_R(s''_m) \zeta'(s''_m,\theta''_m) & \textrm{if } s''_m > +1 \\
\zeta(s''_m,\theta''_m) + \zeta'(s''_m,\theta''_m)
& \textrm{if } -1 \le s''_m \le +1 \\
\gamma_R(-s''_m) \zeta(s''_m,\theta''_m) +
\zeta'(s''_m,\theta''_m) & \textrm{if } s''_m < -1
\end{array} \right.
$$
and $\zeta^0$ coincides with $\zeta$ and $\zeta'$
on the rest of $\Sigma$ and $\Sigma'$ respectively.

Let us define the splitting map
\begin{eqnarray*}
s_R : L_{k-1}^{p,d}(\Lambda^{0,1}(E''))
&\rightarrow& L_{k-1}^{p,d}(\Lambda^{0,1}(E)) \oplus
L_{k-1}^{p,d}(\Lambda^{0,1}(E')) \\
\eta'' &\rightarrow& (\eta,\eta')
\end{eqnarray*}
where, near the punctures $\ux_{\us+1-m}$ and $\ox'_m$,
$$
\left\{ \begin{array}{rcl}
\eta(s_m,\theta_m) &=& \beta(s''_m) \eta''(s''_m,\theta''_m)  \\
\eta'(s'_m,\theta'_m) &=& (1- \beta(s''_m)) \eta''(s''_m,\theta''_m)
\end{array} \right.
$$
and $\eta, \eta'$ coincide with $\eta''$ away from these punctures. \\
Note that the operators $g_R$ and $s_R$ are uniformly bounded in $R$.

Let $Q_\infty$ be a right inverse for the surjective operator
$K \oplus_{\Gamma_t} L$. Let us define an approximate right
inverse $\tilde{Q}_R$ for $M$ using the following commutative diagram :
$$
\begin{CD}
L_{k-1}^{p,d}(\Lambda^{0,1}(E'')) @>\tilde{Q}_R>> \Gamma''_{\os'',\us''}
\oplus L_k^{p,d}(E'') \\
@Vs_RVV    @AAg_RA \\
L_{k-1}^{p,d}(\Lambda^{0,1}(E)) \oplus L_{k-1}^{p,d}(\Lambda^{0,1}(E'))
@>Q_{\infty}>>
\Gamma''_{\os'',\us''} \oplus \Gamma_t \oplus L_k^{p,d}(E) \oplus L_k^{p,d}(E')
\end{CD}
$$
Note that $\tilde{Q}_R$ is uniformly bounded in $R$, since $g_R$ and $s_R$ are.

\begin{proof}[Proof of Proposition~\ref{rinv}]
By construction, $M_R \tilde{Q}_R \eta'' = \eta''$ away from the neighborhoods
of the gluing. On the other hand,
for $s''_m \in [-\frac{R-R_0}2,+\frac{R-R_0}2]$, we have
\begin{eqnarray*}
M_R \tilde{Q}_R \eta'' &=& M_R \zeta'' \\
&=& M_R (v_m + \gamma_R(-s''_m) \zeta + \gamma_R(s''_m) \zeta') \\
&=& M_R v_m + \gamma_R(-s''_m) K \zeta + \gamma_R(s''_m) L \zeta'
- \frac{d}{ds}\gamma_R(-s''_m) \zeta + \frac{d}{ds}\gamma_R(s''_m) \zeta' \\
&& + \gamma_R(-s''_m) (M - K) \zeta + \gamma_R(s''_m) (M - L) \zeta'
\end{eqnarray*}
But
\begin{eqnarray*}
\gamma_R(-s''_m) K \zeta + \gamma_R(s''_m) L \zeta' &=&
\gamma_R(-s''_m) \eta + \gamma_R(s''_m) \eta' \\
&=&  \eta +  \eta' \\
&=& \eta''
\end{eqnarray*}
because $\gamma_R(-s) = 1$ (resp. $\gamma_R(s) = 1$) when $\eta$
(resp. $\eta'$) is not zero.\\
Therefore,
\begin{eqnarray*}
M_R \tilde{Q}_R \eta'' - \eta'' &=& M_R v_m - \gamma'_R(-s''_m) \zeta
+ \gamma'_R(s''_m) \zeta' \\
&& + \gamma_R(-s''_m) (M_R - K) \zeta + \gamma_R(s''_m) (M_R - L) \zeta' .
\end{eqnarray*}
Hence, it follows from our constructions that
$$
\|M_R \tilde{Q}_R \eta'' - \eta''\|_R \le C(R) (|v| + \|\zeta\| + \|\zeta'\|)
$$
where $C(R) \rightarrow 0$ when $R \rightarrow \infty$.

Therefore, if $R$ is sufficiently large,
$$
\| M_R \tilde{Q}_R - I \| \le \frac12 .
$$
Hence, the operator $M_R \tilde{Q}_R$ is invertible.
Let $Q_R = \tilde{Q}_R (M_R \tilde{Q}_R)^{-1}$. By construction,
$Q_R$ is a right inverse for $M_R$, and it is uniformly bounded in $R$.
\end{proof}

\begin{corollary}  \label{gluemap}
There is a natural isomorphism
$$
\phi : \ker K \oplus_{\Gamma_t} \ker L \rightarrow \ker M_R
$$
that is defined up to homotopy, if $R$ is sufficiently large.
\end{corollary}

\begin{proof}
Let $Q_R$ be the uniformly bounded right inverse for the glued operator
$M_R \in {\cal O}''$, as in Proposition~\ref{rinv}.
We define $\phi$ to be the composition of the map $g_R$ and the projection map
$(I - Q_R M_R)$ from the domain of $M_R$ to its kernel, along the image of $Q_R$.
In other words, $\phi = (I - Q_R M_R) \circ g_R$. \\
We claim that the restriction of $\phi$ to $\ker K \oplus_{\Gamma_t}
\ker L$ is an isomorphism for $R$ sufficiently large. By Proposition~\ref{Fred},
the dimensions of both spaces agree, so it is enough to show that
$\phi$ is injective. \\
By contradiction, assume that for any large $R$, we can find $\zeta_R \in
\ker K$ and $\zeta'_R \in \ker L$ such that
$\| \zeta_R \| + \| \zeta'_R \| = 1 $ and
$g_R(\zeta_R,\zeta'_R) = Q_R \eta''_R$, for some
$\eta''_R$. Since $K \zeta_R = L \zeta'_R = 0$, we see that
$\lim_{R \rightarrow \infty} M_R g_R(\zeta_R, \zeta'_R) = 0$ by the same kind of computation
as for $M_R \tilde{Q}_R \eta''$ at the beginning of the proof of Proposition~\ref{rinv}
But since $M_R Q_R = I$, it follows that $\eta''_R \rightarrow 0$
when $R \rightarrow \infty$. Using this in the original equation gives
$\lim_{R \rightarrow \infty} g_R(\zeta_R, \zeta'_R) = 0$. But this contradicts
$\| \zeta_R \| + \| \zeta'_R \| = 1$.
\end{proof}

Suppose now that the operator $K \oplus_{\Gamma_t} L$ is not surjective. Then
we stabilize the operators $K$ and $L$ using finite dimensional oriented vector
spaces $H_1$ and $H_2$ generated by smooth sections with compact support of
$\Lambda^{0,1}(E)$ and $\Lambda^{0,1}(E')$ respectively. We obtain operators
\begin{eqnarray*}
K_{H_1} : H_1 \oplus \Gamma_{\os,\us} \oplus L^{p,d}_k(E)  &\longrightarrow&
L^{p,d}_{k-1}(\Lambda^{0,1}(E)) \\
(h_1,v,\zeta) &\longrightarrow& h_1 + K(v,\zeta) \\
L_{H_2} : H_2 \oplus \Gamma'_{\os',\us'} \oplus L^{p,d}_k(E')  &\longrightarrow&
L^{p,d}_{k-1}(\Lambda^{0,1}(E')) \\
(h_2,v',\zeta') &\longrightarrow& h_2 + L(v',\zeta').
\end{eqnarray*}
It is a standard property of determinant spaces (see for example
\cite{Floer/Hofer:orient}) that $\Det(K)$ and $\Det(L)$ are
canonically isomorphic to $\Det(K_{H_1})$ and $\Det(L_{H_2})$, respectively.
We can always choose $H_1$ and $H_2$ so that
$K_{H_1} \oplus_{\Gamma_t} L_{H_2}$ is surjective. When $R$ is large enough, the
operator $K_{H_1} \oplus_{\Gamma_t} L_{H_2}$ is a stabilization of
$K \oplus_{\Gamma_t} L$; hence, their determinant spaces are canonically isomorphic.

Moreover, the symplectic orientation of the first summand of
$E_{\ux_{\us +1 - m}} \cong E_{\ox'_m} \cong \RR^2 \oplus \RR^{2n-2}$
induces an orientation on the vector space $\Gamma_t$. This in turn induces a
natural isomorphism between $\Det(K \oplus_{\Gamma_t} L)$ and the determinant space
of its stabilization $K \oplus L$.

Therefore, we can generalize Corollary~\ref{gluemap} as follows :

\begin{corollary} \label{glueorient}
Given operators $K$ and $L$ as above, there is a natural isomorphism
$$
\Psi : \Det(K) \otimes \Det(L) \longrightarrow \Det(M_R)
$$
that is defined up to homotopy.
\end{corollary}

\section{Coherent Orientations of the Determinant Line Bundles}
\label{orientation}

We will explain the algorithm to orient all determinant line bundles over the
spaces  ${\cal O}(\Sigma,E;\oA_1,\ldots,\oA_{\os};\uA_1,\ldots,\uA_{\us})$
of Fredholm operator. This was originally done in \cite{Floer/Hofer:orient}
for $\os=\us=1$.

In order to construct the orientations, we need the following operations :
\begin{enumerate}
\item[(i)] Gluing of orientations. If we have elements
\begin{eqnarray*}
K & \in &
{\cal O}(\Sigma,E;\oA_1,\ldots,\oA_{\os};\uA_1,\ldots,\uA_{\us})\\
L & \in & {\cal O}
(\Sigma',E';\oA'_1,\ldots,\oA'_{\os'};\uA'_1,\ldots,\uA'_{\us'})
\end{eqnarray*}
with $\uA_{\us +1-i}$ matching $\oA'_{i}$ for $i=1,\ldots,t\le
\min\{\us,\os' \}$ we can glue them to obtain a new element
$$
K\sharp_t L \in {\cal O}
(\Sigma'',E'';\oA_1,\ldots,\oA_{\os},\oA'_{t+1},\ldots,
\oA'_{\os'};\uA_1,\ldots,\uA_{\us -t},\uA'_1,\ldots,\uA'_{\us'})
$$
with the corresponding glued data. By Corollary~\ref{glueorient},
given orientations $o_K$ and $o_L$ of the determinant bundles of $L$
and $K$ we get an induced orientation $o_K \sharp_t o_L$ of the
determinant of $K \sharp_t L$.  In addition this operation is associative.

\item[(ii)] Disjoint union of orientations.
In addition to gluing we may form the disjoint union $L_1 \amalg L_2$ of
two operators $L_1$ and $L_2$. There is a natural isomorphism from
$\Det(L_1)\otimes\Det(L_2)\to \Det(L_1\amalg L_2)$. We denote by
$o_{L_1}\amalg o_{L_2}$ the induced orientation of $\Det(L_1\amalg L_2)$
This operation may look very trivial but  it is not. In fact, it depends
on the  order of the components $L_1$ and $L_2$.
\end{enumerate}

Now we are ready to set up the necessary orientation algorithm.

{\bf Step 0. }
For each closed Riemann surface $\Sigma$, we can choose a $\CC$-linear
element $L \in {\cal O}(\Sigma,E;\emptyset;\emptyset)$, namely a
$\delbar$--operator on $E$. The kernel and
cokernel of this operator are complex vector spaces, and hence are canonically
oriented. We define $o_L$ to be this complex orientation.

{\bf Step 1. } For each closed Reeb orbit $\gamma \in {\cal P}_{\alpha}$,
there is a unique map $\gamma : S^{1} \longrightarrow M$
parametrizing it with constant arc length such that $\gamma(1)$ maps
to the point $z_{\gamma}$. In addition we fix a symplectic trivialization
of the normal bundle of $\gamma$ (i.e.~of $\xi$ along the orbit).
Thus we obtain a (periodic) family of symplectic matrices
$A_{\gamma}(t)$ by linearizing the Reeb flow.
We pick the trivial symplectic vector bundle $\theta^n$ of rank
$n$ over the Riemann sphere $\CC P^{1\ast}$ with one positive puncture
identifying the restriction of $\overline{\theta}$ to the boundary with
$\xi|_{\gamma}$ using the trivialization.
We fix an element
$\oL_{\gamma} \in {\cal O}(\CC P^{1},\theta;A_{\gamma};\emptyset)$
(note that these sets are all non--empty) together with an orientation
$o_{\oL_\gamma}$ of the determinant of the differential operator
$\oL_{\gamma}$.


{\bf Step 2. } For each closed Reeb orbit $\gamma \in {\cal P}_{\alpha}$
pick the trivial symplectic vector bundle $\theta^n$ of rank $n$
over the Riemann sphere $\CC P^{1\ast}$ with one negative puncture and fix an
element $\uL_{\gamma} \in {\cal O}(\CC P^1,\theta;\emptyset;A_{\gamma})$.
We define the orientation $o_{\uL_\gamma}$ of its determinant so that
$o_{\uL_\gamma} \sharp_1 o_{\oL_\gamma}$ coincides with the complex
orientation induced by the usual $\delbar$--operator on $\CC P^1$.

{\bf Step 3. } For $L_i \in {\cal O}(\Sigma,E;A_i;\emptyset)$, $i = 1, \ldots, k$,
we define $o_{L_1 \coprod \ldots \coprod L_k} = o_{L_1}\amalg \ldots \amalg o_{L_k}$.
We have a similar definition for negative punctures.
For each
$$
K \in {\cal O}(\Sigma,E;\oA_1,\ldots,\oA_{\os};\uA_1,\ldots,\uA_{\us}) ,
$$
we define the orientation $o_K$ of its determinant so that
$$
o_{\uL_{\og_{\os}}  \coprod \ldots \coprod \uL_{\og_1}}
\sharp_{\os} o_{K} \sharp_{\us} o_{\oL_{\ug_{\us}} \coprod \ldots
\coprod \oL_{\ug_1}}
$$
coincides with the complex orientation induced by the usual
$\delbar$--operator on $\Sigma$.

Note that step 3 is compatible with steps 1 and 2, i.e. the orientations
$o_{\oL_\gamma}$ and $o_{\uL_\gamma}$ obtained using step 3, coincide with
the ones defined in step 1 and 2, respectively, for each closed Reeb orbit
$\gamma \in {\cal P}_\alpha$.

Now we have to determine the behavior of these orientations under the gluing
operation $\sharp_t$. Let us first describe the situation of two
surfaces $\Sigma$ and $\Sigma'$ and operators
$K\in {\cal O}(\Sigma,E;\oA_1,\ldots,\oA_{\os}; \uA_1,\ldots,\uA_t)$,
$L\in {\cal O}(\Sigma',E'; \oA'_{t},\ldots,\oA'_1;\uA'_1,
\ldots,\uA'_{\us})$ with
$\uA_m$ matching $\oA'_m$, for $m = 1, \ldots, t$.
That means we are in a situation of a
complete gluing, i.e.~$K \sharp_t L\in {\cal O}(\Sigma'', E'';
\oA_1, \ldots, \oA_{\os}; \uA'_1, \ldots, \uA'_{\us})$.

\begin{proposition} \label{totglue}
In the situation described above we have
$$
o_K \sharp_t o_L = o_{K \sharp_t L}.
$$
\end{proposition}

\begin{proof}
Note that the statement is true by definition if $t=1$, when
$\Sigma=\CC P^{1}$ and $\os =0$ or when $\Sigma'=\CC P^{1}$ and $\us' =0$,
since we simply cap off a puncture.

Next, note that the statement also holds
if $\os =1$ or $\us' =1$ instead of zero, i.e.~in the case of
gluing a cylinder. Indeed, we can reduce this case to the previous one
after capping off the other end of the cylinder, by associativity of the gluing
operation $\sharp$ for orientations.

Consider now the case of Riemann surfaces with an arbitrary topology and arbitrary
number of punctures, equipped with operators $K$ and $L$ such that
their asymptotic expression near the punctures are complex linear.
The operators $K$ and $L$ are then homotopic to a
complex linear operator in the class of operators with
fixed asymptotics at the punctures. Hence, their determinants inherit
a natural orientation from their complex structure.
The same is of course true for the glued
operator and it is clear that gluing complex orientations
yields the complex orientation since the gluing map is homotopic
to a complex linear isomorphism in this case.

Finally, we reduce the general case to the above
situation by gluing cylinders $Z$ and $Z'$ to the punctures at which
we want to glue $K$ and $L$. We choose $Z$ and $Z'$ so that
the glued operators $K \sharp Z$ and $Z'\sharp L$ have complex linear
asymptotic expressions at the negative and positive punctures, respectively.
Since we can cap off all positive or negative ends of $K$ and $L$, we
may assume that there are no such punctures.
We have a homotopy
$K \sharp_t L \cong K \sharp_t Z \sharp_t Z' \sharp_t L$
and hence
$$
o_{K \sharp_t L} = o_{K \sharp_t Z \sharp_t Z' \sharp_t L} = o_{K \sharp_t Z}
\sharp_t o_{Z'\sharp_t L}
$$
due to consistency of gluing in the complex linear case. The latter is
equal to
$$
(o_{K \sharp_t Z}) \sharp_t (o_{Z'}\sharp_t o_L) =
((o_{K \sharp_t Z}) \sharp_t o_{Z'}) \sharp_t o_L
=o_{K \sharp_t Z\sharp_t Z'} \sharp_t o_L = o_K
\sharp_t o_L .
$$
\end{proof}

We are ready to derive the property for reordering the punctures.

\begin{proof}[Proof of Theorem~\ref{orientorderings}] This
property follows from the consistency of the
orientations under disjoint union applied to the $1$--punctured
spheres used to cap off operators over arbitrary surfaces with
punctures. Comparing the orientations via $\overline{\iota}_{l,l+1}$ is the
same as comparing the orientations on the disjoint unions
$o_{\oL_l} \amalg o_{\oL_{l+1}}$ and $o_{\oL_{l+1}} \amalg o_{\oL_l}$
of $\oL_l \in {\cal O}(\CC P^1,\theta;\oA_{\ug_l};\emptyset)$
and $\oL_{l+1} \in {\cal O}(\CC P^1,\theta;\oA_{\ug_{l+1}};\emptyset)$.
This gives an orientation of the (virtual) vector space
$\ker(\oL_l) \oplus \ker(\oL_{l+1}) \ominus (\coker(\oL_l) \oplus
\coker(\oL_{l+1}))$
which is given by orientations of the two vector spaces of the difference.
Thus we can compare the orientations induced by either of the
two possible orders of the components with that given one.
It is not hard to check that this differs exactly if both
indices $\mbox{ind}(\oL_l)$ and $\mbox{ind}(\oL_{l+1})$ are odd. But
$$
\mbox{ind}(\oL_l) = \mu({\ug_l})+(n-1) .
$$
The proof for the case of changing the ordering at the
positive ends is almost the same; in that case, we have to
use spheres with one negative puncture and the index is
$$
\mbox{ind}(\uL_k) = -(\mu({\og_k})+(n-1)) .
$$
\end{proof}

Next we discuss the behavior of coherent orientations under the
disjoint union operation.

\begin{proposition} \label{union}
Let
$$
K \in {\cal O}(\Sigma,E;\oA_1,\ldots,\oA_{\os};\uA_1,\ldots,\uA_{\us})
$$
and
$$
L \in {\cal O}(\Sigma',E';\oA'_1,\ldots,\oA'_{\os'};
\uA'_1,\ldots,\uA'_{\us'}).
$$
Then, for
$$
K \amalg L \in {\cal O}(\Sigma \amalg \Sigma', E \amalg E';
\oA_1, \ldots, \oA_{\os}, \oA'_1, \ldots,
\oA'_{\os'};\uA_1,\ldots,\uA_{\us}, \uA'_1,\ldots,\uA'_{\us'}),
$$
we have $o_{K \coprod L} = (-1)^\epsilon o_L \amalg o_K$, where
$$
\epsilon = \Big( \sum_{l=1}^{\us} (\mu(\ug_l) + (n-3)) \Big)
\Big( \sum_{k=1}^{\os'} (\mu(\og'_k) + (n-3)) \Big).
$$
\end{proposition}

\begin{proof}
By definition,
$$
o_{\uL_{\og'_{\os'}}\coprod\ldots\coprod \uL_{\og'_1} \coprod
\uL_{\og_{\os}}\coprod\ldots\coprod \uL_{\gamma_1}}\sharp_{\os + \os'}
o_{K \coprod L}
\sharp_{\us + \us'} o_{\oL_{\ug_{\us'}} \coprod\ldots\coprod
\oL_{\ug'_1} \coprod \oL_{\ug_{\us}}\coprod\ldots\coprod
\oL_{\ug_1}}
$$
coincides with the complex orientation. On the other hand,
\begin{eqnarray*}
&& o_{\uL_{\og'_{\os'}}\coprod\ldots\coprod \uL_{\og'_1} \coprod
\uL_{\og_{\os}}\coprod\ldots\coprod \uL_{\gamma_1}}\sharp_{\os + \os'}
o_K \amalg o_L
\sharp_{\us + \us'} o_{\oL_{\ug_{\us'}} \coprod\ldots\coprod
\oL_{\ug'_1} \coprod \oL_{\ug_{\us}}\coprod\ldots\coprod
\oL_{\ug_1}} \\
&=&
o_{\uL_{\og'_{\os'}}\coprod\ldots\coprod \uL_{\og'_1}}
\sharp_{\os'}
o_{\uL_{\og_{\os}}\coprod\ldots\coprod \uL_{\og_1} \sharp_{\os} K}
\amalg o_{L \sharp_{\us'} \oL_{\ug'_{\us'}} \coprod\ldots\coprod
\oL_{\ug'_1}} \sharp_{\us}
o_{\oL_{\ug_{\us}}\coprod\ldots\coprod \oL_{\ug_1}} \\
&=&
(-1)^\epsilon
o_{\uL_{\og'_{\os'}}\coprod\ldots\coprod \uL_{\og'_1}}
\sharp_{\os'}
o_{L \sharp_{\us'} \oL_{\ug'_{\us'}} \coprod\ldots\coprod
\oL_{\ug'_1}} \amalg
o_{\uL_{\og_{\os}}\coprod\ldots\coprod \uL_{\og_1}\sharp_{\os} K}
\sharp_{\us} o_{\oL_{\ug_{\us}}\coprod\ldots\coprod \oL_{\ug_1}}
\end{eqnarray*}
But the latter is a complex orientation multiplied by $(-1)^\epsilon$, by
definition of $o_K$ and $o_L$. Comparison with the definition of $o_{K \coprod L}$
gives the desired result.
\end{proof}

The following statement about the behavior of coherent orientations for
a general gluing is a mixed case of Propositions~\ref{totglue}  and
\ref{union}.

\begin{corollary} \label{glue}
Let
$$
K \in {\cal O}(\Sigma,E;\oA_1,\ldots,\oA_{\os};\uA_1,\ldots,\uA_{\us})
$$
and
$$
L \in {\cal O}(\Sigma',E';\oA'_1,\ldots,\oA'_{\os'};
\uA'_1,\ldots,\uA'_{\us'})
$$
with $\uA_{\us +1-m}$ matching $\oA'_m$ for $m = 1,\ldots,t\le
\min\{\us,\os'\}$.
Then, for
$$
K \sharp_t L \in {\cal O}(\Sigma \sharp_t \Sigma', E \sharp_t E';
\oA_1, \ldots, \oA_{\os}, \oA'_{t+1}, \ldots,
\oA'_{\os'};\uA_1,\ldots,\uA_{\us-t},
\uA'_1,\ldots,\uA'_{\us'}) ,
$$
we have
$o_{K \sharp_t L} = (-1)^\epsilon o_K \sharp_t o_L$, where
$$
\epsilon = \Big( \sum_{l=1}^{\us-t} (\mu(\ug_l) + (n-3)) \Big)
\Big( \sum_{k=t+1}^{\os'} (\mu(\og'_k) + (n-3)) \Big) .
$$
\end{corollary}

\begin{proof}
By definition,
$$
o_{\uL_{\og'_{\os'}}\coprod\ldots\coprod \uL_{\og'_{t+1}}
\coprod \uL_{\og_{\os}}\coprod\ldots\coprod \uL_{\og_1}}
\sharp_{\os + \os' -t} o_{K \sharp_t L} \sharp_{\us +\us' -t}
o_{\oL_{\ug'_{\us'}} \coprod\ldots\coprod \oL_{\ug'_1}
\coprod \uL_{\ug_{\us -t}}\coprod\ldots\coprod L'_{\ug_1}}
$$
coincides with the complex orientation.
On the other hand,
$$
o_{\uL_{\og'_{\os'}}\coprod\ldots\coprod \uL_{\og'_{t+1}}
\coprod \uL_{\og_{\os}}\coprod\ldots\coprod \uL_{\og_1}}
\sharp_{\os +\os' -t} (o_K \sharp_t o_L) \sharp_{\us +\us' -t}
o_{\oL_{\ug'_{\us'}} \coprod\ldots\coprod \oL_{\ug'_1}
\coprod \uL_{\ug_{\us -t}}\coprod\ldots\coprod L'_{\ug_1}}
$$
coincides with
$$
(o_{\uL_{\og'_{\os'}}\coprod\ldots\coprod \uL_{\og'_{t+1}}} \amalg
o_{\uL_{\og_{\os}}\coprod\ldots\coprod \uL_{\og_1}
\sharp_{\os} K}) \sharp_{\us + \os' -t}
(o_{L \sharp_{\us'} \oL_{\ug'_{\us'}} \coprod\ldots\coprod
\oL_{\ug'_1}}
\amalg o_{\oL_{\ug_{\us -t}}\coprod\ldots\coprod \oL_{\ug_1}})
$$
because of associativity. In order to apply Proposition~\ref{totglue}
to the gluing $\sharp_{\us +\os'-t}$, we have to permute several punctures in
the right term : $\oA'_{t+1}, \ldots, \oA'_{\os'}$ must be exchanged
with $\uA_{\us -t},\ldots, \uA_1$.
This gives precisely the sign $(-1)^\epsilon$.
After the gluing $\sharp_{\us +\os'-t}$, we obtain a complex orientation, and
comparison with the definition of $o_{L \sharp_t K}$ gives the
desired result.
\end{proof}

Finally, we study the behavior of coherent orientations under the action of
automorphisms.

\begin{proposition} \label{autom}
Let $L \in {\cal
O}(\Sigma,E;\oA_1,\ldots,\oA_{\os};\uA_1,\ldots,\uA_{\us})$ and
$\sigma$ be a diffeomorphism of $\Sigma$ preserving the punctures.
Let $j' = \sigma_* j$ and $K = \sigma(L)$. Then the map induced by
$\sigma$ from the determinant line of $L$ to the determinant line
of $K$ preserves coherent orientations.
\end{proposition}

\begin{proof}
First, if $\Sigma$ is a closed Riemann surface, we can deform $L$ into a
$\CC$-linear operator $L'$ by a path of operators. The
isomorphism induced by $\sigma$ from the kernel and cokernel of $L'$ to the
kernel and cokernel of $\sigma(L')$ is
$\CC$-linear, so it sends $o_{L'}$ to $o_{\sigma(L')}$ and hence $o_L$ to
$o_K$.
In the general case, consider the action of $\sigma$ on
$o_{\uL_{\og_{\os}}\coprod\ldots\coprod \uL_{\og_1}}\sharp_{\os}
o_{L} \sharp_{\us} o_{\oL_{\ug_{\us}}\coprod\ldots\coprod
\oL_{\ug_1}}$. On one hand, $\sigma$ preserves the complex orientation
as above. On the other hand, the action of $\sigma$ clearly commutes with the
gluing maps, and by assumption $\sigma$ extends as the identity on each capping
1-punctured sphere. Therefore the coherent orientations must be preserved as
well.
\end{proof}

This proposition shows that the determinant line bundle for the Fredholm
operators remains trivial if we allow the conformal structure on $\Sigma$
to vary. Moreover, for $j' = j$ and $K = L$, it shows that the coherent
orientations descend to the quotient space, since they are preserved by
the automorphisms.

\begin{remark} \label{anom}
One is tempted to
consider the space of Fredholm operators with asymptotics in a
connected component of ${\cal S}$ instead of concrete descriptions.
Unfortunately, this is not possible, due to Theorem~2 in
\cite{Floer/Hofer:orient}. Neither of the determinant bundles over
such spaces is orientable !
\end{remark}

\begin{remark}
We can compare the coherent orientations constructed in this paper with
the orientations introduced in section 1.8 of
\cite{Eliashberg/Givental/Hofer:SFT}. Even though
our definition of coherent orientations is formulated differently than
axioms (C1)-(C3) of \cite{Eliashberg/Givental/Hofer:SFT}, the 2 notions actually
coincide, provided we choose the $1 \otimes 1^*$ orientation for each cylinder
equipped with a translation invariant operator
$L \in {\cal O}(\CC P^1,\theta;A_\gamma;A_\gamma)$.
\end{remark}

\section{Moduli Spaces of $J$-holomorphic curves}  \label{moduli}

We will apply our construction of coherent orientations to moduli
spaces of holomorphic curves in a symplectic cobordism
$(W^{2n},\omega)$ of dimension $2n$. We describe the Banach
setting in which we study holomorphic curves in  symplectic
manifolds with ends and show that the linearization of the
$\delbar$--equation gives rise to a Fredholm operator between
Banach spaces. The approach we finally chose was explained to us
by Helmut Hofer.

As usual we denote by ${\cal T}_{g,s}$, the Teichm\"uller space of
genus $g$ smooth complex curves with $s$ distinct marked points.
This is a smooth complex manifold of complex dimension $3g-3+s$.
We represent ${\cal T}_{g,s}$ by a family of conformal structures
on a (fixed) closed oriented surface $\Sigma$ and $s$ distinct
points. We therefore denote by $(j,x_1,...,x_s)\in {\cal T}_{g,s}$
the corresponding data on $\Sigma$.

Let $p > 2$ and $k \ge 1$. We denote by
$$
{\cal B}_{W,J}^{\Sigma}
(\og_1,\dots,\og_{\os}; \ug_1,\dots,\ug_{\us})
$$
the set of elements consisting of $(j,\ox_1, \ldots,
\ox_{\os},\ux_1,\ldots, \ux_{\us})\in {\cal T}_{g,\os+\us}$  and a
map
$u:\Sigma\setminus\{\ox_1,\ldots,\ox_{\os},\ux_1,\ldots,\ux_{\us}\}
\longrightarrow W$ which is locally in $L^p_k$. For each puncture
there is a neighborhood which is completely mapped inside one of
the ends of $W$ : for a positive puncture this will be the
positive end, for a negative puncture it will be the negative end.
Furthermore, we ask that the projection of the map to the contact
manifold corresponding to such an end maps a small neighborhood of
$\ox_k$ or $\ux_l$ completely inside a fixed tubular neighborhood
of the closed Reeb orbit $\og_k \in {\cal P}_{\alpha_+}$ or $\ug_l
\in {\cal P}_{\alpha_-}$. We introduce coordinates
$(\vartheta,\zeta) \in S^1 \times D^{2n-2}$ in a tubular
neighborhood of the geometric locus $|\gamma|$ of a closed
$\alpha_\pm$--Reeb orbit $\gamma$ with multiplicity $m(\gamma)$,
such that $\{0\}\times S^1 = |\gamma|$ parameterizes the
underlying simple closed Reeb orbit with constant arc length
element. We require that $\vartheta = 0$ corresponds to the fixed
point on $\gamma$.

Finally, near a puncture we require that there are holomorphic
cylindrical coordinates $(s,\theta)\in [R_0,\infty) \times S^1$
for a positive puncture and $(s,\theta)\in (-\infty,-R_0] \times
S^1$ for a negative puncture, such that the components of the map
$u=(u_\RR, u_M)$, with $u_M=(\vartheta, \zeta)$, satisfy $(u_{\RR}
- {\cal A}(\gamma)s - a_{0}),(\vartheta - m(\gamma)\theta), \zeta
\in L^p_k(e^{d|s|}ds \, d\theta)$ in the tubular coordinates
around the corresponding Reeb orbit.

Here ${\cal A}(\gamma)=\int_\gamma \alpha_\pm$ is the action of
the (multiple) closed Reeb orbit $\gamma$, $a_{0}$ is a constant,
$0 < d/p < 1$ is smaller than the absolute value of the biggest
negative eigenvalues of $L^\infty_{S_{\og_k}}$ for $\ox_k$ or
smaller than the smallest positive eigenvalue of
$L^\infty_{S_{\ug_l}}$ for $\ux_l$, and $L^\infty_{S_{\gamma}}$ is
an operator acting on smooth functions in the following way : fix
a hermitian trivialization of $\xi$ over $\gamma$, i.e.~trivialize
it as a complex and symplectic bundle. Then
$L^\infty_{S_{\gamma}}:= -\mbox{i}\frac{\partial}{\partial
\theta}-S_{\gamma}(\theta)$ where $S_{\gamma}(\theta) =
-\mbox{i}\pi_{m(\gamma)\theta}dR(m(\gamma)\theta)
\pi_{m(\gamma)\theta}$ is a family of symmetric matrices, $\pi :
TM\longrightarrow \xi$ is the orthogonal projection (see
\cite{Hofer/Wysocki/Zehnder:asymptotics}). Obviously, maps in
${\cal B}_{W,J}^{\Sigma} (\og_1,\dots,\og_{\os};
\ug_1,\dots,\ug_{\us})$ are asymptotic to $\og_k$ in $\ox_k$ at
$+\infty$ and to $\ug_l$ in $\ux_l$ at $-\infty$.

Remember that the asymptotics impose a condition on the choice of
complex cordinates near the punctures. Namely, there is a finite
set of $m(\gamma)$ possible {\em directions}
$\frac{\partial}{\partial x}$ at the punctures, since $\lim_{x\to
0}u_M(x)=z_\gamma$ for $x\in \RR$, due to our choices. Finally, we
consider as an element in ${\cal B}$ the map $u$ together with a
choice of possible directions at each puncture. We may think of
the additionally introduced direction at the punctures as a
puncture on the boundary of the compactification of the Riemann
surface to a surface with boundary using the limits.

In \cite{Hofer/Wysocki/Zehnder:asymptotics} the authors show that
punctured holomorphic curves with {\em finite energy} which appear
naturally in SFT are asymptotic to Reeb orbits at the punctures
and satisfy the above decay conditions  such that they lie in the
set ${\cal B}$. Moreover, with its obvious topology this set
becomes a Banach manifold. This can be seen using the exponential
map on $W$, corresponding to a complete Riemannian metric which is
cylindrical on the ends of $W$.

Let ${\cal Hol}_{W,J}^{\Sigma}(\og_1,\dots,\og_{\os},
\ug_1,\dots,\ug_{\os})$ be the set of punctured holomorphic maps
with the given asymptotics. The above discussion shows that ${\cal
Hol}_{W,J}^{\Sigma} \subset {\cal B}_{W,J}^{\Sigma}$, so that
${\cal B}$ is an appropriate configuration space to contain the
set of holomorphic maps. Notice that in order to obtain the moduli
space ${\cal M}_{W,J}^{\Sigma}(\og_1,\dots,\og_{\os},
\ug_1,\dots,\ug_{\os})$ of punctured holomorphic curves we still
need to divide out by the mapping class group $\mbox{Diff}(\Sigma,
x)/\mbox{Diff}_0(\Sigma,x)$.

The tangent space $T_{(j,x,u)}{\cal B}$ splits into the tangent
space of ${\cal T}_{g,\os + \us}$ at $(j,x)$ and the vector space
of sections of $u^{\ast}TW$ which are locally $L^p_k$ and, near
the punctures, are either $L^{p,d}_k(TW)$ in cylindrical
coordinates or are constant linear combinations of
$\frac{\partial}{\partial t}$ and $R_{\alpha_\pm}$.

At each point $(j,x,u) \in {\cal B}$ we have a Banach space
$$
E_{(j,x,u)}:= L^{p,d}_{k-1}(\Lambda^{0,1}(u^{\ast}TW))
$$
forming a smooth Banach bundle ${\cal E}\longrightarrow {\cal B}$.
Then the equation for $J$--holomorphic curves defines
a Fredholm section $\delbar:{\cal B}\longrightarrow {\cal E}$.

Consider the linearization $\delbar_{(j,x,u)}$ of this section at
$(j,x,u) \in {\cal B}$. This operator splits into
$$
\delbar_{(j,x,u)} : D_{(j,x,u)} \oplus L_{(j,x,u)} :
T_{(j,x)}{\cal T}_{g,\os+\us} \oplus \Gamma_{\os,\us} \oplus
L^{p,d}_k(u^\ast TW) \longrightarrow
L^{p,d}_{k-1}(\Lambda^{0,1}(u^\ast TW)) .
$$
It follows from  the above discussion that the second summand $L_{(j,x,u)}$ is an
element of ${\cal O}(\Sigma,E;\oA_1,\ldots,\oA_{\os};\uA_1,\ldots,\uA_{\us})$.

Therefore, by Proposition~\ref{Fred}, $L_{(j,x,u)}$ is Fredholm and has index
$$
\mbox{ind}(L_{(j,x,u)}) = \sum_{k=1}^{\os}(\mu(\og_{k})-(n-1))
- \sum_{l=1}^{\us}(\mu(\ug_{l})+(n-1)) + n\chi(\Sigma)
+ 2c_{1}(W,\omega)[A].
$$

Since the first summand $D_{(j,x,u)}$ of $\delbar_u$ has finite rank,
from the index above one can easily derive the expected dimension of the
moduli spaces :
$$
\sum_{k=1}^{\os} \mu(\og_{k})
- \sum_{l=1}^{\us} \mu(\ug_{l}) + (n-3)(\chi(\Sigma)-\os-\us)
+ 2c_{1}(W,\omega)[A].
$$

\begin{proof}[Proof of Theorem~\ref{gluing}]
First, we have to explain how to orient the determinant bundles of
the moduli spaces $\cal{M}_{W,J}^{A,\Sigma}(\og_1,\dots,\og_{\os},
\ug_1,\dots,\ug_{\os})$ using the construction of the previous
section.

Let $(j,x,u) \in {\cal Hol}_{W,J}^{A,\Sigma}(\og_1,\dots,\og_{\os},
\ug_1,\dots,\ug_{\us})$. Consider the linearized operator
$\delbar_{(j,x,u)}$ splitting as $D_{(j,x,u)} \oplus L_{(j,x,u)}$. The family of
operators $L_{(j,x,u)}$ for $(j,x,u) \in {\cal B}$ defines a continuous map
\begin{eqnarray*}
op : {\cal Hol}_{W,J}^{A,\Sigma}(\og_1,\dots,\og_{\os},
\ug_1,\dots,\ug_{\us}) &\longrightarrow& {\cal O}
(\Sigma, u^*TW; A_{\og_1},\dots,A_{\og_{\os}};
A_{\ug_1},\dots,A_{\ug_{\us}}) \\
(j,x,u) &\longrightarrow& L_{(j,x,u)} .
\end{eqnarray*}

The operator $\delbar_{(j,x,u)}$ is homotopic to $0\oplus
L_{(j,x,u)}$, a stabilization of the operator $L_{(j,x,u)} \in
{\cal O}$ with the complex vector space $T_{(j,x)}{\cal
T}_{g,\os+\us}$. Hence, their determinant spaces are canonically
isomorphic and the determinant bundle over ${\cal Hol}$
corresponding to the Fredholm operator $\delbar_{(j,x,u)}$ is
isomorphic to $op^* \Det({\cal O})$. In particular, the
orientations from Section \ref{orientation} induce orientations on
the determinant bundle of ${\cal Hol}$.

The mapping class group of $(\Sigma,x)$ acts naturally on ${\cal
Hol}$. Let ${\cal M}$ be the moduli space of this action. By
Proposition~\ref{autom}, the coherent orientations are preserved
by this action and descend to the moduli spaces ${\cal M}$.

Next, we have to show that, given compact subsets
$$
K_1 \subset {\cal M}^{\Sigma_1}_{W_1,J_1} (\og_1,\dots,\og_{\os};
\ug_1,\dots,\ug_{\us-t},\beta_1,\dots, \beta_t)
$$
and
$$
K_2 \subset {\cal M}^{\Sigma_2}_{W_2,J_2}
(\beta_t,\dots,\beta_2,\beta_1,\og'_{t+1}, \dots,\og'_{\os'};
\ug'_1,\dots,\ug'_{\us'}) ,
$$
the gluing maps
$$
\Phi_R : K_1 \times K_2 \longrightarrow
{\cal M}^{\Sigma}_{W_R,J_R}
(\og_1,\dots,\og_{\os},\og'_{t+1},\dots,\og'_{\os'};
\ug_1,\dots,\ug_{\us-t},\ug'_1, \dots,\ug'_{\us'})
$$
are orientation preserving up to a sign $(-1)^\epsilon$ where
$\epsilon$ is determined by Corollary~\ref{glue}.

In order to prove this, we have to relate the differential of $\Phi_R$ to the
linear gluing map of Corollary~\ref{gluemap}. Note that the construction of the
gluing map requires the Fredholm section $\delbar : {\cal B} \rightarrow
{\cal E}$ to be transverse to the zero section. In other words, the
linearization $\delbar_{(j,x,u)}$ of this section at $(j,x,u) \in
\delbar^{-1}(0) \subset {\cal B}$ must be a surjective operator.
This transversality can be achieved in different ways, such as perturbing the
almost complex structure $J$ or perturbing the right hand side of the
Cauchy-Riemann equation with an element of
$L^{p,d}_{k-1}(\Lambda^{0,1}(u^{\ast}TW))$ having compact support in the
complementary of the punctures. In all cases, the linearization of the perturbed
section $\delbar$  is still an element of ${\cal O}$ modulo some zero order
terms not affecting our constructions. Therefore, the definition of coherent
orientations is independent of
the precise way we achieve transversality for the moduli spaces.

Let $(j_i,x_i,u_i) \in K_i$ for $i = 1, 2$. Using cutoff functions
near the closed orbits $\beta_1, \ldots, \beta_t$, we construct a
pre-glued map $u_R$ into the glued cobordism $(W_R,J_R)$. Gluing
the two conformal structures we also construct  conformal data
$(j_R,x_R)$ on the glued surface $\Sigma_1 \sharp_t \Sigma_2$.
Then  the glued data $(j_R,x_R,u_R)$ satisfy
$$
\| \delbar u_R \|\le C(R) ,
$$
where $C(R)\to 0$ as $R\to\infty$.

Under the above transversality assumption, the linearized operator
$\delbar_{(j_R,x_R,u_R)}$ is surjective and has a uniformly bounded right
inverse, as in Proposition \ref{rinv}. An actual holomorphic curve
$(j_R,x_R,w_R)$
is obtained using Newton iterations (see \cite{McDuff/Salamon:Quantum}), in a
neighborhood of $(j_R,x_R,u_R)$ of size controlled by $C(R)$. In particular,
the difference in norm of the linearizations $\delbar_{(j_R,x_R,w_R)}$ and
$\delbar_{(j_R,x_R,u_R)}$, and the glued operator $\delbar_R$ obtained from
$\delbar_{(j_i,x_i,u_i)}$, for $i = 1, 2$ as in Section \ref{setup},
are arbitrarily small for $R$ sufficiently large.

The differential of the gluing map $\Phi_R$ is the composition of
the linearizations of the pre-gluing map and of the Newton
iteration map. The linearization of the pre-gluing map involves
gluing sections of $L^{p,d}_k(u_i^* TW)$ for $i = 1, 2$ using some
cut-off functions. It is therefore of the same form as the linear
gluing map $g_R$ from Section \ref{setup}. On the other hand, the
differential of the Newton iteration map approaches the projection
$I - Q_{u_R} \delbar_{u_R}$ to the kernel of $\delbar_{u_R}$ along
the image of its right inverse $Q_{u_R}$ as $R$ becomes large. The
above discussion shows that this projection is very close to the
projection $I - Q_R \delbar_R$ of Section \ref{setup} for $R$
large.

This shows that the differential of $\Phi_R$ and the linear gluing map are very
close for $R$ large, so that they induce the same map on orientations.
\end{proof}

Note that multiple covers of holomorphic curves have to be
considered in Symplectic Field Theory. These prevent, even for a
generic choice of $J$, the linearized Cauchy-Riemann operator from
being surjective. Therefore, it is necessary to use more
sophisticated methods, such as multi-valued perturbations
\cite{Liu/Tian}, in order to achieve transversality. The moduli
spaces we obtain in this way are {\em branched labeled
pseudo--manifolds}. These are topological spaces consisting of
finitely many smooth open strata of finite dimension such that the
top--dimensional strata are dense. There is an assignment of
positive rational numbers to each top--stratum. This discussion
shows that, to define Symplectic Field Theory for all contact
manifolds, we have to leave the realm of integer invariants, and
rather use rational coefficients.

As soon as we wish to compute integral or rational invariants we have to
define orientations  of the moduli spaces of holomorphic curves
which satisfy some coherence property under the gluing operation.
This will lead to signs for holomorphic curves of index $1$ in
symplectizations (i.e.~those used to define the differential)
and for holomorphic curves of index $0$ or $-1$ in symplectic
cobordisms (i.e.~those which define homomorphisms and homotopies
between the differential algebras). Counting them with their signs
we get integers or rational numbers.

On the other hand,
there could be situations where one is able to define ${\ZZ}_{2}$--contact
homology (see for example \cite{Ustilovsky}), and hence avoid coherent orientations.
However, it is not obvious that one can justify the invariance even in the
absence of holomorphic curves of index $1$ : multiply covered cylinders occur
quite naturally in the symplectic cobordism defining the chain map which should be
chain homotopic to the identity.
In that case, multi--valued perturbations could not be avoided, and
rational coefficients would make coherent orientations necessary.

\section{Even and odd behavior of the orientation for multiply-covered
orbits}\label{behavior}

We now establish the behavior of the orientations
under change of asymptotic directions at punctures which are
asymptotic to Reeb orbits with multiplicity.

\begin{proof}[Proof of Theorem~\ref{orientbadorbits}]
First note that, in view of steps 2 and 3 in Section \ref{orientation} for the
construction of coherent orientations, it is enough to prove the theorem for
the moduli space ${\cal M}^{\Sigma}(\gamma_m,\emptyset)$, where $\gamma_m$ is
the $m$-fold covering of the simple orbit $\gamma$.

The index of the operator $\oL_{\gamma}$ is given by $\mu(\gamma) + (n-1)$.
Then $\oL_{\gamma_m}$ is the pullback of $\oL_{\gamma}$ under the
branched covering $z \mapsto z^m$. In particular, this operator has a
$\ZZ_m$--symmetry of rotations. Hence, $\ZZ_m$ acts on the kernel and
cokernel of $\oL_{\gamma_m}$, and these vector spaces split into irreducible
representations over $\RR$.

If $m$ is odd, there are no elements of even order, so the action of $\ZZ_m$
preserves orientations. Hence, coherent orientations are invariant under a
change of asymptotic direction.

If $m$ is even, the possible irreducible representations include 2
representations of dimension 1 (trivial and sign change) and $m-2$
representations of dimension 2, generated by rotations of angle
$\frac{2\pi}m i$, $i = 1, \ldots, m-1$, $i \neq \frac{m}2$. The
trivial representations correspond to the kernel and cokernel of
$\oL_\gamma$. Therefore, the index difference
$$
ind(\oL_{\gamma_m})-ind(\oL_\gamma)=\mu(\gamma_m) -\mu(\gamma)
$$
has the same parity as the multiplicity of the orientation
reversing representation over $\RR$ in the kernel and cokernel of
$\oL_{\gamma_m}$. Hence, rotating the asymptotic direction of
$\frac{2\pi}m$ reverses the coherent orientation if and only if
$\mu(\gamma_m) - \mu(\gamma)$ is odd. This is exactly the case
described in the theorem (see e.g.~Lemma~3.2.4. in
\cite{Ustilovsky:thesis}).
\end{proof}

Coherent orientations can be constructed after we make some orientation choices
for all closed Reeb orbits. We now explain how to reduce the number of
necessary choices in the case of multiply covered orbits.

\begin{proposition} \label{redchoice}
Let $\gamma$ be a simple closed Reeb orbit. Fix an
orientation of the determinant bundle over ${\cal O}
(\CC P^{1},\theta^n;A_{\gamma},\emptyset)$.
\begin{enumerate}
\item[(i)]  If $\gamma'$ is an odd multiple of $\gamma$,
there is a canonical choice for the orientation of the determinant bundle over
${\cal O}(\CC P^1,\theta^n;A_{\gamma'},\emptyset)$
induced by the orientation for $\gamma$.
\item[(ii)] The choice of an orientation for the double orbit $\gamma'$, in
addition to $\gamma$, induces orientations for all multiple orbits.
\item[(iii)] If $n=2$, the orientation for $\gamma$ alone induces
orientations for even multiples as well, if they are good.
\end{enumerate}
\end{proposition}

\begin{proof}
$(i)$ Let $m$ be the odd multiplicity of $\gamma'$. Recall that
$\ZZ_m$ acts by rotation on the kernel and cokernel of
$\oL_{\gamma'}$ and that these vector spaces split into the
invariant part and some irreducible representations over $\RR^2$
generated by a rotation of angle $\frac{2\pi}m i$, $i = 1, \ldots,
m-1$. The invariant part consists of the pullback of the kernel
and cokernel of the operator $\oL_\gamma$ and hence it is
oriented. We choose the orientation on each irreducible
representation $\RR^2$ in such a way that $i < \frac{m}2$.

$(ii)$ If the multiplicity $m$ of $\gamma'$ is even, the kernel
and cokernel of $\oL_{\gamma'}$ split into the invariant part, some irreducible
representations over $\RR^2$ generated by a rotation of angle
$\frac{2\pi}m i$, $i = 1, \ldots, m-1$, $i \neq \frac{m}2$, and
the part $V_{-1}$ on which the generator of $\ZZ_m$ acts by $-1$.
The invariant part and the $\RR^2$ summands are oriented as above.
On the other hand, the kernel and cokernel of the operator corresponding to
the double orbit split precisely into the invariant part and $V_{-1}$.
Since the choice of an orientation for the simple orbit is equivalent
to an orientation of the invariant part, the choice of an orientation for the
double orbit is just what we need to orient all the other multiples.

$(iii)$
In the case $n=2$ we have to find a way to orient the part $V_{-1}$
on which the rotation acts by $-1$. Choose the operator $\oL_\gamma$ so
that it splits into the standard Cauchy-Riemann operator $\delbar$ and an
operator $\delbar_{S_\gamma}$ of the form
$$
\frac{\partial}{\partial s} + \mbox{i} \frac{\partial}{\partial t} +
S_{\gamma}(t)
$$
near the positive puncture. Since $\delbar_{S_\gamma}$ has nondegenerate
asymptotics, its Fredholm index is given by $\mu(S_\gamma) - 1$.
This kind of operators were well studied in
\cite{Hofer/Wysocki/Zehnder:Embeddednes}. We now summarize
some facts about $\delbar_{S_\gamma}$ from that paper, assuming
$\textrm{ind}(\delbar_{S_\gamma}) \ge 0$. Otherwise, we consider its formal
adjoint, exchanging kernel and cokernel.

\begin{enumerate}
\item[a)] Let $\phi_\lambda$ be an eigenvector with eigenvalue $\lambda\in \RR$
of the operator
$L^\infty_{S_\gamma}:=\mbox{i}\frac{\partial}{\partial t} +
S_\gamma(t)$.
 It has a well-defined winding number $w(\phi_\lambda)$ in the
fixed trivialization. Each winding number appears exactly twice
and increases with the eigenvalue. Two eigenvectors with the same
winding number are {\em pointwise}  linearly independent in
$\RR^2$ or a multiple of each other.

\item[b)] The operator $\delbar_{S_\gamma}$ is surjective. Elements in
$\ker \delbar_{S_\gamma}$ have the form
$Ce^{s\lambda}\phi_{\lambda} + O(e^{s\lambda})$, where $C\neq 0$,
$\lambda<0$ and $1 \le w(\phi_\lambda) \le w(S_\gamma)$, where
$$
w(S_\gamma)  =  \max\{w(\phi_\lambda)\mid 0>\lambda\in
\mbox{spec}(\mbox{i}\frac{\partial}{\partial t} + S_\gamma(t))\}.
$$

\item[c)]Consersely,  we may choose a basis of $\ker \delbar_{S_\gamma}$
such that the leading terms are in $1-1$ correspondence  to the
set of eigenvectors $\{\phi_\lambda\}_\lambda$ satisfying the
above conditions.

Therefore, the Maslov index of $S_\gamma$ can be computed in terms
of these data. We have
$$
\mu(S_\gamma) = 2w(S_\gamma) + p(S_\gamma)
$$
where $p(S_\gamma)  =  0$ if the winding of the smallest positive eigenvalue is
the same as that of the biggest negative, and $p(S_\gamma)=1$ otherwise.

\end{enumerate}

Let $\gamma'$ be the double cover of the orbit $\gamma$. Let $\oL_{\gamma'}$
be the pullback of the operator $\oL_\gamma$ under the branched covering
$z \mapsto z^2$. The kernel of $\oL_{\gamma'}$ splits into the pullback of
$\ker \oL_\gamma$ and the part $V_{-1}$ that is not invariant under the
$\ZZ_2$ action. By assumption, $\gamma'$ is good, so that $p(S_{\gamma'}) =
p(S_\gamma)$. Therefore, $V_{-1}$ splits into summands $\RR^2$ generated by
2 eigenvectors of equal winding number. Since these 2 eigenvectors are
pointwise linearly independent, the corresponding summand $\RR^2$ inherits a
natural orientation from the complex orientation of $\xi$.
Hence, we obtain a natural orientation on $V_{-1}$.
\end{proof}

\begin{remark}
If we reduce the number of choices in the construction of coherent
orientations using Proposition \ref{redchoice}, we can extract
more precise information from the moduli spaces. Consider the
symplectic cobordism $(\RR \times M, d(e^t \alpha_s))$, where
$\alpha_s$ interpolates between contact forms $\alpha_0$ and
$\alpha_1$ for isomorphic contact structures $\xi_0 \cong \xi_1$.
Counting holomorphic curves of index $0$ in this cobordism, we
obtain an isomorphism $\Phi$ between the contact homologies
computed using $\alpha_0$ and $\alpha_1$. Assume the image of a
generator $\gamma$ has the form $\Phi(\gamma) = c
\widetilde{\gamma}$, where $c \neq 0$. If we choose orientations
separately for each closed orbit, then of course the sign of $c$
is meaningless. However, if we choose orientations using
Proposition \ref{redchoice}, then the sign of $c$ may contain some
extra information.
\end{remark}

\end{document}